\documentclass[12pt,reqno]{amsart}
\usepackage[headings]{fullpage}
\usepackage{amssymb,amsmath,mathtools,bbm,tikz,
tikz-cd,color,relsize,multirow}
\usepackage{makecell}
\usepackage[all,cmtip]{xy}
\usepackage{url}
\usepackage[group-separator={,},group-minimum-digits={3}]{siunitx}
\usetikzlibrary{positioning}
\usetikzlibrary{calc}
\usetikzlibrary{decorations.markings}
\usetikzlibrary{arrows}
\usetikzlibrary{calc}
\usetikzlibrary{decorations.markings}
\tikzstyle{hvector}=[inner sep=2pt,draw=blue!50,fill=blue!10,thick]
\tikzstyle{unit}=[inner sep=2pt,shape=circle, draw]
\tikzstyle{counit}=[inner sep=2pt,shape=circle, draw,fill=gray]
\tikzstyle{antipode}=[inner sep=2pt,shape=rectangle, draw]
\tikzstyle{cocycle}=[inner sep=2pt,shape=circle, draw]
\tikzstyle{twistedm}=[inner sep=2pt,shape=circle, fill=gray]
\tikzstyle{autom}=[inner sep=2pt,shape=circle, draw]
\tikzstyle{coact}=[inner sep=2pt,shape=circle, fill=black]

\usepackage[bookmarks=true,%
    colorlinks=true,%
    linkcolor=blue,%
    citecolor=blue,%
    filecolor=blue,%
    menucolor=blue,%
    urlcolor=blue,%
    breaklinks=true]{hyperref}
\usepackage{slashed}    
\usepackage{listings}        
\usepackage{verbatim}
\usepackage[normalem]{ulem}  
\usepackage{diagbox}
\usepackage{booktabs}
\usepackage[shortlabels]{enumitem}

\theoremstyle{definition}

\def\be{\begin{equation}}
\def\ee{\end{equation}}

\newcommand{\HTW}{\mathrm{HTW}}

\definecolor{codegreen}{rgb}{0,0.6,0}
\definecolor{codegray}{rgb}{0.5,0.5,0.5}
\definecolor{codepurple}{rgb}{0.58,0,0.82}
\definecolor{backcolour}{rgb}{0.95,0.95,0.92}

\lstdefinelanguage{PARIGP}{
  keywords={typeof, new, true, false, catch, function, return, null, catch,
    switch, var, if, in, while, do, else, case, break},
  keywordstyle=\color{blue}\bfseries,
  ndkeywords={class, export, boolean, throw, implements, import, this},
  ndkeywordstyle=\color{darkgray}\bfseries,
  identifierstyle=\color{black},
  sensitive=false,
  comment=[l]{//},
  morecomment=[s]{/*}{*/},
  commentstyle=\color{purple}\ttfamily,
  stringstyle=\color{red}\ttfamily,
  morestring=[b]',
  morestring=[b]"
}

\lstdefinestyle{code}{
    backgroundcolor=\color{backcolour},   
    commentstyle=\color{codegreen},
    keywordstyle=\color{magenta},
    numberstyle=\tiny\color{codegray},
    stringstyle=\color{codepurple},
    basicstyle=\ttfamily\footnotesize,
    breakatwhitespace=false,         
    breaklines=true,                 
    captionpos=b,                    
    keepspaces=true,                 
    numbers=left,                    
    numbersep=5pt,                  
    showspaces=false,                
    showstringspaces=false,
    showtabs=false,                  
    tabsize=2
}

\begin{document}

\title[Table of mutant knots with at most 16 crossings]{
Table of mutant knots with at most 16 crossings}
\author{Stavros Garoufalidis}
\address[S. Garoufalidis]{
  International Center for Mathematics, Department of Mathematics \\
  Southern University of Science and Technology \\
  Shenzhen, China \newline
  {\tt \url{http://people.mpim-bonn.mpg.de/stavros}}
  }
\email{stavros@mpim-bonn.mpg.de}

\author{Shana Yunsheng Li}
\address[S.Y. Li]{
  Department of Mathematics \\
  University of Illinois \\
  Urbana, IL, USA \newline
  {\tt \url{https://shana-y-li.github.io}}\newline
  ORCID: {\tt \url{https://orcid.org/0009-0006-5320-106X}}
  }
\email{yl202@illinois.edu}

\author{Daniel Tubbenhauer}
\address[D. Tubbenhauer]{
  School of Mathematics and Statistics \\
  The University of Sydney \\
  F07 Carslaw Building, Office 827 \\
  NSW 2006, Australia \newline
  {\tt  \url{https://www.dtubbenhauer.com}} \newline
  ORCID: {\tt \url{https://orcid.org/0000-0001-7265-5047}}
  }
\email{daniel.tubbenhauer@sydney.edu.au}

\author{Victor L. Zhang}
\address[V.L. Zhang]{
  School of Mathematics and Statistics \\
  University of New South Wales (UNSW) \\
  NSW 2052, Australia \newline
  {\tt \url{https://dustbringer.github.io}} \newline
  ORCID: {\tt \url{https://orcid.org/0009-0007-5799-6477}}
  }
\email{victor.l.zhang@unsw.edu.au}

\subjclass[2020]{Primary 57K10, 68T09; Secondary 57K14, 57M12.}

\thanks{
  {\em Key words and phrases:}
  knots, Conway mutation, knot tables.
}

\date{\today}

\begin{abstract}
We give a complete, duplication-free table of all Conway mutant
equivalence classes of prime knots (and their mirrors) up to 16 
crossings along
with certificates of mutation and non-mutation.
\end{abstract}

\maketitle

{\footnotesize
\tableofcontents
}


\section{Introduction}
\label{sec.intro}

Knot tabulation has a long and fascinating history told
time and again; see for example Przytycki's presentation 
\cite{Przytycki}.

In the nineteen seventies, Conway introduced mutation~\cite{Conway}
of knots, closely related to knot tabulation, and to 
skein-theory, predating both hyperbolic 
geometry~\cite{Thurston,HTW} and quantum topology
~\cite{RT:ribbon}. Among Conway's discoveries was the mutant
pair $(11n34,11n42)$ of 11-crossing knots with trivial 
Alexander polynomial. 

Conway mutation was then comprehensively studied 
by Bonahon--Siebenmann~\cite{Bonahon} using characteristic 
submanifolds, 
analogue of the JSJ decomposition of 3-manifolds. A table of Conway 
mutant prime knots with at most 15 crossings was later given by
Stoimenow using a mixture of topological invariants
~\cite{St:hard, St:tabulating}. 

Our aim is 
to compute the Conway mutation classes of all prime knots with
at most 16 crossings, distinguishing chiral knots from their
mirrors. We give explicit mutations connect the members of each 
class,
and mutation invariants separate distinct classes. This extends
Stoimenow's census through 15 crossings~\cite{St:tabulating}.

To phrase our results, 
let $\HTW$ denote the set of \num{1701935} prime knots 
with at most 16 crossings \cite{HTW}. Adding their mirrors and
removing duplications results into a set $\HTW^+$ of
\num{3401978} knots. The map $\HTW^+ \to \HTW$
is onto and either $1:1$ or $2:1$ depending on 
a knot being amphicheiral or chiral. 
We partition $\HTW$ and $\HTW^+$ into Conway equivalence classes 
$\langle \HTW \rangle$ and $\langle \HTW^+ \rangle$.
So, we have a
commutative diagram
\begin{equation*}
\begin{tikzcd}
\HTW^+ \arrow[r, two heads] \arrow[d] & \HTW 
\arrow[d] \\
\langle \HTW^+ \rangle \arrow[r, two heads] & \langle \HTW \rangle
\end{tikzcd}
\end{equation*}
where the bottom map is onto and either $1:1$ or $2:1$. 
Our mutation table and certificates are available in 
\cite{knotmutations}, which contains both the plain data files and a 
python module \texttt{knot\_mutations} as an add-on to 
\texttt{Spherogram} and \texttt{SnapPy}. \autoref{tab:3} summarizes 
our results. 

\begin{table}[htbp]
\centering
\small
\begin{tabular}{@{}c r r@{}}
\toprule
Class size & $\langle\HTW\rangle$ & $\langle\HTW^+\rangle$ \\
\midrule
 1 & \num{1216217} & \num{2430792} \\
 2 &   \num{171268} &   \num{342522} \\
 3 &    \num{6012} &    11,824 \\
 4 &    25,677 &    \num{51454} \\
 5 &         0 &         0 \\
 6 &     1,483 &     2,962 \\
 7 &        10 &        20 \\
 8 &     1,481 &     2,962 \\
 9 &         2 &         6 \\
10 &         8 &         8 \\
11 &         0 &         0 \\
12 &       111 &       222 \\
13-15 &      0 &         0 \\
16 &        12 &        28 \\
\bottomrule
\end{tabular}
\vspace{10pt}
\caption{Number of Conway mutation equivalence classes
of each size in $\HTW$ and $\HTW^+$.}
\label{tab:3}
\end{table}

Along the way, we found an omitted Conway mutant
pair $(15a13571,15a13629)$ in Stoime\-nov's original table. 
(This was concurrently observed by Werner Alois Stanggassinger, 
assisted by AI.) 

In our data, there are 130 chiral knots that are mutant to their 
mirrors, and all their mutation classes include at least one 
amphicheiral knot. There is no example of two knots with different 
crossing numbers that are mutant in our data.

\smallskip
\noindent \textbf{Acknowledgments.} The authors wish to thank 
Nathan Dunfield for enlightening conversations.
The authors thank the Keeling cluster maintenance team at UIUC, 
the Katana cluster 
(DOI: \url{https://doi.org/10.26190/669X-A286}) of 
Research Technology Services at UNSW Sydney and the cluster 
in the Zelmanov Mathematics Center in SUSTech  
for providing computational resources.
SG and DT thank Hans Boden, Zsuzsanna Dancso, and
Iva Halacheva for organizing the Dror-Fest
\textit{A Panorama of Quantum Topology} in
July 2026, where this collaboration began. DT thanks the 
organizers of that workshop for bringing them to a place where 
the mountains were considerably easier to distinguish than 
the knots.

\smallskip
\noindent \textbf{AI declaration.}
Large language models were used for 
writing code and translating human-written code between languages.

\section{Methodology}
\label{sec.method}

We start with the table of all \num{1701935} prime knots with at 
most 16 crossings \cite{HTW}. Adding the mirrors of these knots,
and removing the duplicate amphicheiral knots, we got 
an initial set $\HTW^+$ of \num{3401978} knots. 

We remark that the chirality of hyperbolic knots is decided by  
\texttt{SnapPy} and that non-hyperbolic knots in the knot 
table up to 16 crossings are chiral (indeed, for torus knots 
this is well-known, see e.g., \cite{Murasugi:textbook} 
and for satellites it is due to \cite{HTW}). 

To partition the set $\HTW^+$
into possible-mutation equivalence classes, we use several known 
Conway mutation invariants 
(the references give mutation invariance) from: 
\begin{itemize}
\item 
classical topology (such as the 2-fold branched cover $\Sigma_2(K)$ 
of a knot $K$~\cite[Prop.3.8.2]{Kawauchi}) and hence
the homology of its $k$-fold covers, or epimorphisms of
its fundamental group to finite simple groups, computed by
\texttt{SageMath} \cite{sagemath}, \texttt{SnapPy} \cite{snappy},
    \texttt{GAP} \cite{GAP4} and \texttt{Magma} \cite{magma},
\item 
hyperbolic geometry (such as the volume~\cite{Ruberman}), 
computed by \texttt{SnapPy}, and 
\item 
quantum topology (such as the colored Jones 
polynomial~\cite{MortonTraczyk}, the
Links--Gould polynomial~\cite{DG,GL:doesnot}, odd Khovanov 
homology over $\mathbb{Z}$~\cite{Bloom}, and
reduced even Khovanov homology over
$\mathbb{Q}$~\cite{KWZ:kh}), computed
by \cite{Li:FPT-RT, GL:patterns}, \texttt{SnapPy} and 
\texttt{KnotJob} \cite{knotjob}.
\end{itemize}
Using these invariants, we partition the set $\HTW^+$ into
equivalence classes of possibly-mutant knots. We remark that we 
also used the fact that if $\Sigma_2(K)$ is hyperbolic, then 
there is no essential Conway sphere in the exterior of $K$, 
hence $K$ is not mutant to any other knot. We thank Nathan 
Dunfield for providing code that certifies the 
hyperbolicity of $\Sigma_2(K)$.

On the other side, we search for mutation certificates to confirm 
that the knots in each of our equivalence classes are mutant to each 
other. This is done with the following programs. The programs and 
scripts can be found at \cite{knotmutations-certificates}.
\begin{itemize}
\item A \texttt{Rust} program that enumerates all possible 
    mutations on a given planar diagram (PD) presentation of a knot.
\item We use \texttt{Regina}'s \cite{regina} PD rewrite function 
    to increase the number of crossings from the starting PD.
\item A \texttt{Rust} program that extends \texttt{Regina}'s 
rewrite function to distinguish mirrored PDs.
\item We use \texttt{SnapPy} \cite{snappy} to identify the knots 
that appear in one given rewrite.
\item For alternating knots we used Greene's method \cite{Greene}, 
which is a complete invariant on the mutation classes of alternating 
knots.
\end{itemize}
The output is a collection of certificates: each certificate 
indicates the mutation between a pair of knots. 

We compute the invariants and mutation certificates described 
above until every pair of knots is known to (respectively) 
be either not mutant or mutant.

\bibliographystyle{hamsalpha}
\bibliography{biblio}
\end{document}